\documentclass[10pt]{article}

\usepackage{amssymb,latexsym}
\usepackage{epsfig}
\usepackage{eufrak}
\usepackage{amsmath}
\usepackage{mathrsfs}
\usepackage{color}

\usepackage{float}

\newtheorem{theorem}{Theorem}
\newtheorem{lemma}{Lemma}
\newtheorem{corollary}{Corollary}

\newcommand{\be}{\begin{equation}}
\newcommand{\ee}{\end{equation}}
\newcommand{\bee}{\begin{eqnarray*}}
\newcommand{\eee}{\end{eqnarray*}}
\newcommand{\bel}{\begin{eqnarray}}
\newcommand{\eel}{\end{eqnarray}}
\newcommand{\bec}{\begin{cases}}
\newcommand{\eec}{\end{cases}}
\newcommand{\bem}{\begin{bmatrix}}
\newcommand{\eem}{\end{bmatrix}}

\newcommand{\la}{\label}
\newcommand{\li}{\left}
\newcommand{\ri}{\right}

\newcommand{\lf}{\lfloor}
\newcommand{\rf}{\rfloor}

\newcommand{\lm}{\lambda}

\newcommand{\de}{\delta}

\newcommand{\al}{\alpha}
\newcommand{\ba}{\beta}

\newcommand{\f}{\frac}
\newcommand{\sq}{\sqrt}
\newcommand{\cd}{\cdots}

\newcommand{\qu}{\quad}
\newcommand{\qqu}{\qquad}

\newcommand{\mscr}{\mathscr}

\newcommand{\bb}{\mathbb}

\newcommand{\wh}{\widehat}
\newcommand{\wt}{\widetilde}
\newcommand{\mrm}{\mathrm}
\newcommand{\bs}{\boldsymbol}

\newcommand{\ap}{\approx}

\newcommand{\sh}{\slash}

\newcommand{\tx}{\text}

\newcommand{\pa}{\partial}

\newcommand{\bed}{\begin{description}}
\newcommand{\eed}{\end{description}}
\newcommand{\bei}{\begin{itemize}}
\newcommand{\eei}{\end{itemize}}
\newcommand{\ben}{\begin{enumerate}}
\newcommand{\een}{\end{enumerate}}
\newcommand{\bib}{\bibitem}
\newcommand{\beL}{\begin{lemma}}
\newcommand{\eeL}{\end{lemma}}
\newcommand{\beT}{\begin{theorem}}
\newcommand{\eeT}{\end{theorem}}
\newcommand{\beC}{\begin{corollary}}
\newcommand{\eeC}{\end{corollary}}

\newcommand{\bpf}{\begin{pf}}
\newcommand{\epf}{\end{pf}}
\newcommand{\bsk}{\bigskip}

\newcommand{\pfbox}{\hfill\mbox{$\Box$}}

\newenvironment{pf}{\paragraph*{Proof{\rm.}}}{\pfbox\bigskip}

\begin{document}

\title{{\bf Probability Estimation with Truncated Inverse Binomial Sampling} }


\author{ Xinjia Chen\\
Department of Engineering Technology\\
Northwestern State University, Natchitoches, LA 71497\\
Tel: (318)357-5521  \quad Fax: (318) 357-6145}

\date{  }

\maketitle

\begin{abstract}

In this paper, we develop a general theory of truncated inverse
binomial sampling. In this theory, the fixed-size sampling and
inverse binomial sampling are accommodated as special cases. In
particular, the classical Chernoff-Hoeffding bound is an immediate
consequence of the theory.  Moreover, we propose a rigorous and
efficient method for probability estimation, which is an adaptive
Monte Carlo estimation method based on truncated inverse binomial
sampling. Our proposed method of probability estimation can be
orders of magnitude more efficient as compared to existing methods
in literatures and widely used softwares.

\end{abstract}


\section{Introduction}

In science and engineering, it is an ubiquitous problem to estimate
the probability of event based on Monte Carlo simulation.  For
instance, in engineering technology, a critical concern is the
probability of failure or risk, which is generally considered as the
probability that certain pre-specified requirements for the relevant
system are violated in the presence of uncertainties. Ever since the
advent of modern computers, extensive research works have been
devoted to quantitative approaches of risk evaluation for
engineering systems (see, e.g., \cite{Gul, DT2006, CRC, Fishman,
Legay, Mitzenmacher, Tep} and the references therein). In additional
to theoretical development, many softwares have been developed for
risk evaluation.  For example, for control systems, a software
called RACT has been developed for evaluating the risk of uncertain
systems \cite{IEEE, IFAC}. Many softwares such as  APMC
\cite{HAPMC}, PRISM \cite{Kwia}, UPPAAL \cite{Alexandre}, have been
developed for evaluating the risk of stochastic discrete event
systems (see, \cite{Gul} and the references therein).

One of the remarkable achievements of existing theories and
softwares is the rigorous control of error in the estimation of
probability, that is, the probability of relevant event can be
evaluated with certified reliability. Theoretically, for a priori
given $\al, \; \de \in (0, 1)$, existing methods are able to produce
an estimate $\wh{p}$ for the true value of the probability $p$ so
that one can be $100 (1 - \de)\%$ confident that $| \wh{p} - p | <
\al$ holds. Unfortunately, existing methods suffer from huge
computational complexity as the margin of absolute error $\al$ is
small, e.g. $10^{-3}$.  The sobering fact is that in many
applications, the probability of failure can be of orders of
magnitude smaller than $10^{-3}$. In order to manage such an
extremely small but critical probability of failure, we must
evaluate it within an error much smaller than its value. Since the
evaluation of probability of failure is the basis for refinement of
an engineering design to achieve a lower risk, it is of paramount
mount to improve its efficiency without scarifying rigorousness.

The remainder of the paper is organized as follows.  In Section 2,
we review existing methods of probability estimation.
 In Section 3, we first develop a general theory of truncated inverse binomial sampling.   After, we
 derive error control theory of fixed size sampling from the general
 theory.  Finally, we
 derive error control theory for inverse binomial sampling from the general
 theory.  In Section 4, we
 propose a new method of probability estimation. In Section 5, we
 investigate the worst-case efficiency of the proposed method of
 probability estimation.  Section 6 is the conclusion.

In this paper, we use the following notations.  The set of positive
integers is denoted by $\bb{N}$.   The probability of event is
denoted by $\Pr \{ . \}$.  The mathematical expectation is denoted
by $\bb{E} [. ]$.  The other notations will be made clear as we
proceed.

\section{Conventional Methods for Probability Estimation}

In most engineering applications, the relevant event that the
specified requirements are not satisfied is so complicated that it
is impossible to calculate its probability in any deterministic
sense.   As the advent of modern computers, the Monte Carlo
simulation techniques have been developed for the estimation of
probability.   The general procedure is as follows:

(i) Build a mathematical model for the event with uncertainties
represented by random variables.

(ii) Use a computer program to implement the mathematical model.

(iii) Generate independent and identically distributed samples for
the random variables  and check the occurrence of the random event
for each sample.

Each generation of random sample and checking of the occurrence of
event is considered as one simulation.  For $n = 1, 2, \cd$, let
$S_n$ denote the number of occurrences among $n$ simulations. The
ratio of $S_n$ to $n$  is called the relative frequency.

According to the law of large numbers, the relative frequency,
\[
\wh{p}_n = \f{S_n}{n}, \]
 tends to the probability $p$ as the
number of simulations $n$ tends to infinity. To control the accuracy
in the evaluation of the probability, an absolute error criterion is
frequently used. Specifically, for a priori given margin of absolute
error $\al \in (0, 1)$ and confidence parameter $\de \in (0, 1)$,
one seeks a positive number $n$ such that, with a probability higher
than $1 - \de$, the relative frequency $\wh{p}_n$ differs from the
true value $p$ by an amount at most $\al$, that is,
 \be \la{statement}
 \Pr \li \{   \li |  \wh{p}_n - p \ri |
< \al \ri \} > 1 - \de \ee in a mathematic language.

In order to satisfy (\ref{statement}), the number of simulations $n$
has to be sufficiently large.  For the sake of efficiency, $n$ is
expected to be as small as possible.   This gives rise to the
critical question: how large $n$ should be to ensure
(\ref{statement}) ?

An approximate answer to this question is derived from the central
limit theorem \cite{Desu}, which suggests that the number of
simulations should be chosen as \be \la{apformula}
 n \ap \f{  Z^2  } { 4 \al^2},
\ee where $Z$ is the positive number such that the integration of
the function $\f{1}{\sq{ 2 \pi} } e^{- x^2 \sh 2}$ from $Z$ to
infinity is equal to $\f{\de}{2}$.   The major criticism of the
formula (\ref{apformula}) is that one never knows the difference
between the desired confidence level $1 - \de$ and the probability
of $\li | \wh{p}_n - p \ri | < \al$ (see, e.g., \cite{Fishman,
Hampel} and the references therein).

The best rigorous result in the literature is due to Chernoff and
Hoeffding \cite{Chernoff, Hoeffding}, which asserts that
$(\ref{statement})$ holds provide that \be \la{CHbound} n > \f{ \ln
\f{2}{\de} }{ 2 \al^2 }. \ee  It has been well-known that the
Chernoff-Hoeffding bound is the tightest possible to guarantee the
absolute criterion $(\ref{statement})$ when no information of $p$ is
available, which is common in applications. Since the
Chernoff-Hoeffding bound is the tightest and rigorous, it has become
the gold standard for risk evaluation, documented in engineering
handbooks, encyclopedia and textbooks, and research articles.
Specially, it has been applied to control of uncertain systems
\cite{ETa, ETb, DT2006, CRC, Tep}, machine learning \cite{Mitchell,
VPb, VP,  Vidyasagar}, operation research \cite{Fishman},
statistical model checking \cite{Legay, Younes}, randomized
algorithms for probability and computing \cite{Mitzenmacher,
Motwani}, etc. It has been implemented at the heart of the
simulation engines of many softwares \cite{Gul, Alexandre, IEEE,
HAPMC, Kwia, IFAC}.  Although the Chernoff-Hoeffding bound is shown
to be tight,  a huge computational complexity can be induced in
applications.   For example, if $\al = 10^{-6}$ and $\de = 10^{-3}$,
the number of simulations suggested by the Chernoff-Hoeffding bound
(\ref{CHbound}) is $3,800,451,229,772$.  Table 1 at below shows the
required number of simulations for $\de = 10^{-3}$ and various
values of $\al$, calculated with the Chernoff-Hoeffding bound.

{\small
\begin{table}[h]
 \caption{Number of Simulations $N_{\mrm{CH}}$ from Chernoff-Hoeffding Bound  ($\de = 10^{-3}$) } \label{table_mle}
\begin{center}
\begin{tabular}{|c||c|c|c|c|c|}
\hline $\al$ & $10^{-3}$ & $10^{-4}$ & $10^{-5}$ & $10^{-6}$ & $10^{-7}$\\
\hline $N_{\mrm{CH}}$ & $3,800,452$ & $380,045,123$ & $38,004,512,298$ & $3,800,451,229,772$ & $380,045,122,977,105$\\
 \hline
\end{tabular}
\end{center}
\end{table}
}

From Table 1, it can be seen that the required number of simulations
becomes an astronautical number for $\al < 10^{-3}$.  This indicates
that the Chernoff-Hoeffding bound is not applicable to the
evaluation of probability of failure for many engineering
applications.

\section{Truncated Inverse Binomial Sampling} \la{EBP}

In this section, we shall develop a general theory of  truncated
inverse binomial sampling.  The focus of the theory is on the error
control for the estimation of the underlying binomial parameter.

\subsection{General Theory}

The following result is a restatement of Theorem 4.1 of Chen
\cite{Chen}.

\beT

\la{themmainA}

Let $0 < \de < 1$ and $0 < \al < \ba$ with $\f{\al}{\ba} +
\f{\al}{2} \leq \f{1}{2}$. Let $X_1, X_2, \cd$ be a sequence of
independent and identically distributed Bernoulli random variables
having the same distribution as $X$ such that $\Pr \{ X = 1 \} = 1 -
\Pr \{ X = 0 \} = p \in (0, 1)$. Define \[ A = \f{ \ba \ln
\f{2}{\de} }{ \al (1 + \ba) \ln (1 + \ba) + (\ba - \al - \al \ba)
\ln \li ( 1 - \f{\al \ba} {\ba - \al} \ri ) }, \qqu \qqu \qqu B =
\li ( \f{\al}{\ba} + \al \ri ) A.  \]  Let $L$ and $W$ be real
numbers such that $L \geq A$ and $W \geq B$.  Define
\[ \bs{m} = \min \li \{ n \in \bb{N} : n > L \; \tx{or} \; \sum_{i =
1}^n X_i
>  W \ri \}, \qqu \qqu  \qqu \wh{\bs{p}} = \f{1}{\bs{m}} \sum_{i=1}^{\bs{m}} X_i. \] Then,
\[ \Pr \li  \{ | \wh{\bs{p}} - p | < \al \; \; \tx{or} \; \; \li | \f{
\wh{\bs{p}}  - p }{ p } \ri | < \ba \ri \} > 1 - \de.
\]

\eeT

Making use of Theorem \ref{themmainA}, we have derived the following
result.

\beT

\la{themmainB}

Let $0 < \de < 1$ and $0 < \al < \ba$ with $\f{\al}{\ba} +
\f{\al}{2} \leq \f{1}{2}$. Let $X_1, X_2, \cd$ be a sequence of
independent and identically distributed Bernoulli random variables
having the same distribution as $X$ such that $\Pr \{ X = 1 \} = 1 -
\Pr \{ X = 0 \} = p \in (0, 1)$. Define \[ A = \f{\ba}{ (1 + \ba)
\ln (1 + \ba) - \ba} \; \f{ \ln \f{2}{\de} }{ \al  }, \qqu \qqu \qqu
B = \li ( \f{\al}{\ba} + \al \ri ) A.  \]  Let $L$ and $W$ be real
numbers such that $L \geq A$ and $W \geq B$.  Define
\[ \bs{m} = \min \li \{ n \in \bb{N} : n > L \; \tx{or} \; \sum_{i =
1}^n X_i
>  W \ri \}, \qqu \qqu  \qqu \wh{\bs{p}} = \f{1}{\bs{m}} \sum_{i=1}^{\bs{m}} X_i. \] Then,
\[ \Pr \li  \{ | \wh{\bs{p}} - p | < \al \; \; \tx{or} \; \; \li | \f{
\wh{\bs{p}}  - p }{ p } \ri | < \ba \ri \} > 1 - \de.
\]

\eeT

\bpf

To show the theorem, it suffices to show that \be \la{critical88}
\f{ \ba \ln \f{2}{\de} }{ \al (1 + \ba) \ln (1 + \ba) + (\ba - \al -
\al \ba) \ln \li ( 1 - \f{\al \ba} {\ba - \al} \ri ) } \leq \f{\ba}{
(1 + \ba) \ln (1 + \ba) - \ba} \; \f{ \ln \f{2}{\de} }{ \al  }. \ee
This amounts to show
\[
\al (1 + \ba) \ln (1 + \ba) + (\ba - \al - \al \ba) \ln \li ( 1 -
\f{\al \ba} {\ba - \al} \ri ) \geq \al [ (1 + \ba) \ln (1 + \ba) -
\ba ]
\]
or equivalently
\[
(\ba - \al - \al \ba) \ln \li ( 1 - \f{\al \ba} {\ba - \al} \ri )
\geq - \al \ba,
\]
which can be written as
\[
\ln \li ( 1 - x \ri ) \geq  \f{ - x }{  1 -x  }
\]
with $x = \f{\al \ba} {\ba - \al}$.  This inequality holds because
$0 < x = \f{\al \ba} {\ba - \al} < 1$ as a consequence of the
assumption on $\al $ and $\ba$.

\epf

Making use of Theorem \ref{themmainB}, we have derived the following
result.

\beT

\la{themmainBB}

Let $0 < \de < 1$ and $0 < \al < \ba < 1$ with $\f{\al}{\ba} +
\f{\al}{2} \leq \f{1}{2}$. Let $X_1, X_2, \cd$ be a sequence of
independent and identically distributed Bernoulli random variables
having the same distribution as $X$ such that $\Pr \{ X = 1 \} = 1 -
\Pr \{ X = 0 \} = p \in (0, 1)$. Define \[ A = \f{1}{ \ln 4 - 1 } \;
\f{ \ln \f{2}{\de} }{ \al \ba }, \qqu \qqu \qqu B = \li (
\f{\al}{\ba} + \al \ri ) A.  \]  Let $L$ and $W$ be real numbers
such that $L \geq A$ and $W \geq B$.  Define
\[ \bs{m} = \min \li \{ n \in \bb{N} : n > L \; \tx{or} \; \sum_{i =
1}^n X_i
>  W \ri \}, \qqu \qqu  \qqu \wh{\bs{p}} = \f{1}{\bs{m}} \sum_{i=1}^{\bs{m}} X_i. \] Then,
\[ \Pr \li  \{ | \wh{\bs{p}} - p | < \al \; \; \tx{or} \; \; \li | \f{
\wh{\bs{p}}  - p }{ p } \ri | < \ba \ri \} > 1 - \de.
\]

\eeT

\bpf

Define $C(\ba) = \f{\ba^2}{ (1 + \ba) \ln (1 + \ba) - \ba}$ for $\ba
> 0$. We claim that \be
\la{crt998833}
 C(\ba) < C(1) = \f{1}{ \ln 4 - 1 } \qu \tx{for $0 <
\ba < 1$.} \ee
 To show the claim, it suffices to show that
$C(\ba)$ is an increasing function of $\ba > 0$.  For this purpose,
define $\varphi(\ba) = \f{1}{C(\ba)}$ and $h (\ba) = \ba \ln (1 +
\ba) - 2 (1 + \ba) \ln (1 + \ba) + 2 \ba $ for $\ba > 0$. Then,
\[
\varphi^\prime (\ba) = \f{ h (\ba) }{\ba^3} < 0
\]
if $h (\ba) < 0$,  which indeed holds because $h^\prime (\ba) =
\f{\ba}{1 + \ba} - \ln (1 + \ba) < 0 $ for $\ba > 0$.  Therefore,
the claim is proven.  The proof of the theorem can be completed by
applying the established claim and Theorem  \ref{themmainB}.

\epf

\subsection{Fixed-Size Sampling}

In this section, we shall demonstrate that the general theory of
truncated inverse binomial sampling can be applied to derive error
control method for estimating binomial parameter based on fixed-size
sampling.

\subsubsection{Mixed Criterion}

 To estimate the parameter of  Bernoulli distribution with a mixed
 error criterion based on fixed-size sampling, we have the following
 result.

\beT

\la{themmainC}

Let $0 < \de < 1$ and $0 < \al < \ba$ with $\f{\al}{\ba} +
\f{\al}{2} \leq \f{1}{2}$. Let $X_1, X_2, \cd$ be a sequence of
independent and identically distributed Bernoulli random variables
having the same distribution as $X$ such that $\Pr \{ X = 1 \} = 1 -
\Pr \{ X = 0 \} = p \in (0, 1)$.   Define
\[  \wh{p}_n = \f{1}{n} \sum_{i=1}^{n} X_i \]
for $n \in \bb{N}$. Then,
\[ \Pr \li  \{ | \wh{p}_n - p | < \al \; \; \tx{or} \; \; \li | \f{
\wh{p}_n  - p }{ p } \ri | < \ba \ri \} > 1 - \de
\]
provided that
\[
n > \f{ \ba \ln \f{2}{\de} }{ \al (1 + \ba) \ln (1 + \ba) + (\ba -
\al - \al \ba) \ln \li ( 1 - \f{\al \ba} {\ba - \al} \ri ) }.
\]

\eeT

\bpf Consider the context of Theorem \ref{themmainA}. Let $n$ be a
positive integer greater than \[ A = \f{ \ba \ln \f{2}{\de} }{ \al
(1 + \ba) \ln (1 + \ba) + (\ba - \al - \al \ba) \ln \li ( 1 - \f{\al
\ba} {\ba - \al} \ri ) }. \]  Define
\[ L = \f{1}{2} \times \li ( n + \max \{ n - 1, A \} \ri ).
\]
Then, $L > A$ and $\lf L \rf + 1 = n$.  Let $W > L + 1$.  Then,
\[
\bs{m} = n, \qqu \wh{\bs{p}} = \wh{p}_{n}
\]
with
\[
n = \lf L \rf + 1 > A.
\]
It follows from Theorem \ref{themmainA} that \[ \Pr \li  \{ |
\wh{p}_n - p | < \al \; \; \tx{or} \; \; \li | \f{ \wh{p}_n  - p }{
p } \ri | < \ba \ri \} = \Pr \li  \{ | \wh{\bs{p}} - p | < \al \; \;
\tx{or} \; \; \li | \f{ \wh{\bs{p}}  - p }{ p } \ri | < \ba \ri \}
> 1 - \de.
\]

\epf

As a consequence of (\ref{critical88}) and Theorem \ref{themmainC},
we have the following result.

\beT

\la{themmainC889}

Let $0 < \de < 1$ and $0 < \al < \ba$ with $\f{\al}{\ba} +
\f{\al}{2} \leq \f{1}{2}$. Let $X_1, X_2, \cd$ be a sequence of
independent and identically distributed Bernoulli random variables
having the same distribution as $X$ such that $\Pr \{ X = 1 \} = 1 -
\Pr \{ X = 0 \} = p \in (0, 1)$.   Define
\[  \wh{p}_n = \f{1}{n} \sum_{i=1}^{n} X_i \]
for $n \in \bb{N}$. Then,
\[ \Pr \li  \{ | \wh{p}_n - p | < \al \; \; \tx{or} \; \; \li | \f{
\wh{p}_n  - p }{ p } \ri | < \ba \ri \} > 1 - \de
\]
provided that
\[
n > \f{\ba}{ (1 + \ba) \ln (1 + \ba) - \ba} \times \f{  \ln
\f{2}{\de}  }{\al}.
\]

\eeT

As a consequence of (\ref{crt998833}) and Theorem
\ref{themmainC889}, we have the following result.

\beT

\la{themmainCC}

Let $0 < \de < 1$ and $0 < \al < \ba < 1$ with $\f{\al}{\ba} +
\f{\al}{2} \leq \f{1}{2}$. Let $X_1, X_2, \cd$ be a sequence of
independent and identically distributed Bernoulli random variables
having the same distribution as $X$ such that $\Pr \{ X = 1 \} = 1 -
\Pr \{ X = 0 \} = p \in (0, 1)$.  Define
\[  \wh{p}_n = \f{1}{n} \sum_{i=1}^{n} X_i \] for $n \in \bb{N}$. Then,
\[ \Pr \li  \{ | \wh{p}_n - p | < \al \; \; \tx{or} \; \; \li | \f{
\wh{p}_n  - p }{ p } \ri | < \ba \ri \} > 1 - \de
\]
provided that
\[
n > \f{1}{ \ln 4 - 1 } \; \f{ \ln \f{2}{\de} }{ \al \ba }.
\]

\eeT

\subsubsection{Reducing to Chernoff-Hoeffding Bound}

Formally, the famous Chernoff-Hoeffding bound can be stated as
follows.

\beT

\la{themmainCH}

Let $0 < \de < 1$ and $0 < \al < 1$. Let $X_1, X_2, \cd$ be a
sequence of independent and identically distributed Bernoulli random
variables having the same distribution as $X$ such that $\Pr \{ X =
1 \} = 1 - \Pr \{ X = 0 \} = p \in (0, 1)$.   Define
\[  \wh{p}_n = \f{1}{n} \sum_{i=1}^{n} X_i \]
for $n \in \bb{N}$. Then,
\[ \Pr \li  \{ | \wh{p}_n - p | < \al \ri \} > 1 - \de
\]
provided that
\[
n > \f{ \ln \f{2}{\de} }{ 2 \al^2 }.
\]

\eeT

Making use of Theorem \ref{themmainC}, we can readily derive the
famous Chernoff-Hoeffding bound.  Our argument is as follows.

Define \[ \mscr{H}  (u, v) = u \ln \f{u}{v} + ( 1 - u) \ln \f{1 -
u}{ 1 - v}
\]
for $u \in (0, 1)$ and $v \in (0, 1)$.  For $p \in (0, 1)$, let
\[
\ba > \max \li \{  \f{\al}{1 - \al}, \; \f{\al}{p}  \ri \}.
\]
Note that
\[
\Pr \li  \{ | \wh{p}_n - p | < \al \; \; \tx{or} \; \; \li | \f{
\wh{p}_n  - p }{ p } \ri | < \ba \ri \} = \Pr \li  \{ | \wh{p}_n - p
| < \al \ri \}
\]
as a consequence of $\ba > \f{\al}{p}$.  Hence, according to Theorem
\ref{themmainC},  \[ \Pr \li \{ | \wh{p}_n - p | < \al \ri \} = \Pr
\li  \{ | \wh{p}_n - p | < \al \; \; \tx{or} \; \; \li | \f{
\wh{p}_n  - p }{ p } \ri | < \ba \ri \} > 1 - \de
\]
for $\ba > \max \li \{  \f{\al}{1 - \al}, \; \f{\al}{p}  \ri \}$
with $\f{\al}{\ba} + \f{\al}{2} \leq \f{1}{2}$ provided that
\[
n > \f{ \ln \f{2}{\de} }{ \mscr{H} \li (  \f{\al}{\ba} + \al,
\f{\al}{\ba} \ri ) } = \f{ \ba \ln \f{2}{\de} }{ \al (1 + \ba) \ln
(1 + \ba) + (\ba - \al - \al \ba) \ln \li ( 1 - \f{\al \ba} {\ba -
\al} \ri ) }.
\]
As a consequence of $\ba > \max \li \{  \f{\al}{1 - \al}, \;
\f{\al}{p}  \ri \}$, we have
\[
0 < \f{\al}{\ba} + \al < 1.
\]
Since $\mscr{H}  (\lm, \lm) = 0, \; \f{\pa  \mscr{H}  (\lm + \al,
\lm) }{ \pa \al} = 0$ for $\al = 0$, and
\[
\f{\pa^2  \mscr{H}  (\lm + \al, \lm) }{ \pa \al^2} = \f{1}{ (\lm +
\al) (1 - \lm - \al) } \geq 4
\]
for $0 < \al < 1 - \lm$, it follows from Taylor's series expansion
formula  that
\[
\mscr{H} \li (  \f{\al}{\ba} + \al, \f{\al}{\ba} \ri )  \geq 2
\al^2.
\]
Hence, \[ \Pr \li  \{ | \wh{p}_n - p | < \al \ri \}  > 1 - \de \]
provided that
\[ n
> \f{ \ln \f{2}{\de}  }{ 2 \al^2  }.
\]
This is the classical Chernoff-Hoeffding bound.

\subsection{Inverse Binomial Sampling}

In this section, we shall demonstrate that the general theory of
truncated inverse binomial sampling reduces to the theory of inverse
binomial sampling by letting $\al \to 0$.

To estimate the parameter of  Bernoulli distribution with a relative
 error criterion based on inverse binomial sampling, we have the following
 result.

\beT

\la{themmainD}

Let $0 < \de < 1$ and $\ba  > 0$. Let $X_1, X_2, \cd$ be a sequence
of independent and identically distributed Bernoulli random
variables having the same distribution as $X$ such that $\Pr \{ X =
1 \} = 1 - \Pr \{ X = 0 \} = p \in (0, 1)$.  Define
\[ N = \min \li \{ n \in \bb{N} : \sum_{i =
1}^n X_i
>  \f{ (1 + \ba) \ln \f{2}{\de}
}{ (1 + \ba) \ln (1 + \ba)  - \ba } \ri \}, \qqu \qqu \qqu \wt{p} =
\f{1}{N} \sum_{i=1}^{N} X_i. \] Then,
\[ \Pr \li  \{ \li | \f{
\wt{p}  - p }{ p } \ri | < \ba \ri \} \geq 1 - \de.
\]

\eeT

\bpf

Consider the context of Theorem \ref{themmainA}.  Define \[ Z^* = \bec 1 & \tx{if} \;  \li | \f{ \wt{p}  - p }{ p } \ri | < \ba \\
0 &  \tx{otherwise} \eec
\]
and
\[
Z_\al = \bec 1 & \tx{if} \; | \wh{\bs{p}} - p | < \al \; \; \tx{or}
\;
\; \li | \f{ \wh{\bs{p}}  - p }{ p } \ri | < \ba \\
0 &  \tx{otherwise} \eec
\]
for $\al > 0$. Observing that \bee B & = & \li ( \f{\al}{\ba} + \al
\ri ) \times \f{ \ba \ln \f{2}{\de} }{ \al (1 + \ba) \ln (1 + \ba) +
(\ba - \al - \al
\ba) \ln \li ( 1 - \f{\al \ba} {\ba - \al} \ri ) }\\
& = & \f{ (1 + \ba) \ln \f{2}{\de} }{ (1 + \ba) \ln (1 + \ba) +
(\f{\ba}{\al} - 1 - \ba) \ln \li ( 1 - \f{\al \ba} {\ba - \al} \ri )
}  \eee  and that \[ \li ( \f{\ba}{\al} - 1 - \ba \ri ) \ln \li ( 1
- \f{\al \ba} {\ba - \al} \ri ) \to - \ba
\]
as $\al \to 0$, we have  \[ B \to \f{ (1 + \ba) \ln \f{2}{\de} }{ (1
+ \ba) \ln (1 + \ba)  - \ba }
\]
as $\al \to 0$.  Hence,
\[
W \to \f{ (1 + \ba) \ln \f{2}{\de} }{ (1 + \ba) \ln (1 + \ba)  - \ba
}
\]
as $\al \to 0$.  It follows that $\bs{m} \to N$ almost surely as
$\al \to 0$.  By the convergence of $\bs{m}$ and the definition of
$Z_\al$, we have that $Z_\al \to Z^*$ almost surely as $\al \to 0$.
By Theorem \ref{themmainA},
\[
\lim_{\al \downarrow 0} \bb{E} [ Z_\al ] = \lim_{\al \downarrow 0}
\Pr \li  \{ | \wh{\bs{p}} - p | < \al \; \; \tx{or} \; \; \li | \f{
\wh{\bs{p}}  - p }{ p } \ri | < \ba \ri \} \geq 1- \de.
\]
 Note that $Z_\al$ is a bounded
random variable.  Applying the bounded convergence theorem (or
monotone convergence theorem), we have
\[
\Pr \li  \{ \li | \f{ \wt{p}  - p }{ p } \ri | < \ba \ri \} = \bb{E}
[ Z^* ] = \lim_{\al \downarrow 0} \bb{E} [ Z_\al ] = \lim_{\al
\downarrow 0} \Pr \li  \{ | \wh{\bs{p}} - p | < \al \; \; \tx{or} \;
\; \li | \f{ \wh{\bs{p}}  - p }{ p } \ri | < \ba \ri \} \geq 1- \de.
\]
Hence, we have the desired result.

\epf

\section{New Method of Probability Estimation}

In this section, we shall propose a novel method for probability
estimation with substantially lower computational complexity, while
guaranteeing rigorousness.

\subsection{A Mixed Criterion}

We propose to use a mixed criterion for probability estimation based
on the following philosophy:

In applications, the accuracy of probability estimation can be
measured in terms of absolute error or relative error, while the
margin of absolute error can be substantially smaller than the
margin of relative error.

Specifically, for a priori given margin of absolute error $\al \in
(0, 1)$, margin of relative error $\ba \in (0, 1)$ and confidence
parameter $\de \in (0, 1)$, we seek a procedure which produces an
estimate $\wh{\bs{p}}$ such that with a probability higher than $1 -
\de$, either the absolute error $| \wh{\bs{p}} - p |$ is less than
$\al$ or the relative error {\small $\li | \f{ \wh{\bs{p}} - p }{ p
} \ri |$} is less than $\ba$, that is, \be \la{mixedc}
 \Pr \li \{   \li | \wh{\bs{p}} - p \ri | < \al  \qu
\tx{or} \qu \li | \f{ \wh{\bs{p}} - p  }{ p }  \ri | < \ba  \ri \}
> 1 - \de.
\ee The mixed criterion is a relaxation of the absolute error
criterion. It is important to note that when $p < \f{\al}{\ba}$, the
mixed criterion is equivalent to the absolute criterion.

\subsection{Rectangular Random Walk}

Making use of Theorem \ref{themmainB}, we define a rectangular
domain and a random walk starting from its left-lower vertex and
through such domain as shown by Figure 1.  The slope of the line
connecting the exit point and the staring vertex is taken as an
estimate of the risk. Specifically, define
\[
 L  = \f{\ba}{ (1 + \ba)
\ln (1 + \ba) - \ba} \; \f{ \ln \f{2}{\de} }{ \al  }, \qqu \qqu W =
\li ( \f{\al}{\ba} + \al \ri ) L.
\]
As in Section 2, for $n \in \bb{N}$, let $S_n$ be the number of
occurrence of event among $n$ simulations.  Let $\bb{B}$ be a
rectangular domain in the $(x, y)$-plane of length $L$ and width $W$
which consists of points with coordinate $(x, y)$ satisfying $0 \leq
x \leq L$ and $0 \leq y \leq W$, that is, \[ \bb{B} = \{ (x, y): 0
\leq x \leq L, \; 0 \leq y \leq W \}.
\]
The stopping rule is as follows:

{\it Observe the tuple $(n, S_n)$ for $n = 1, 2, \cd$ until $n$
reach some number $\bs{m}$ such that $(\bs{m}, S_{\bs{m}})$ falls
outside of the box $\bb{B}$, that is, $\bs{m} > L$ or $S_{\bs{m}} >
W$. }

\begin{figure}[H]
\centerline{\psfig{figure=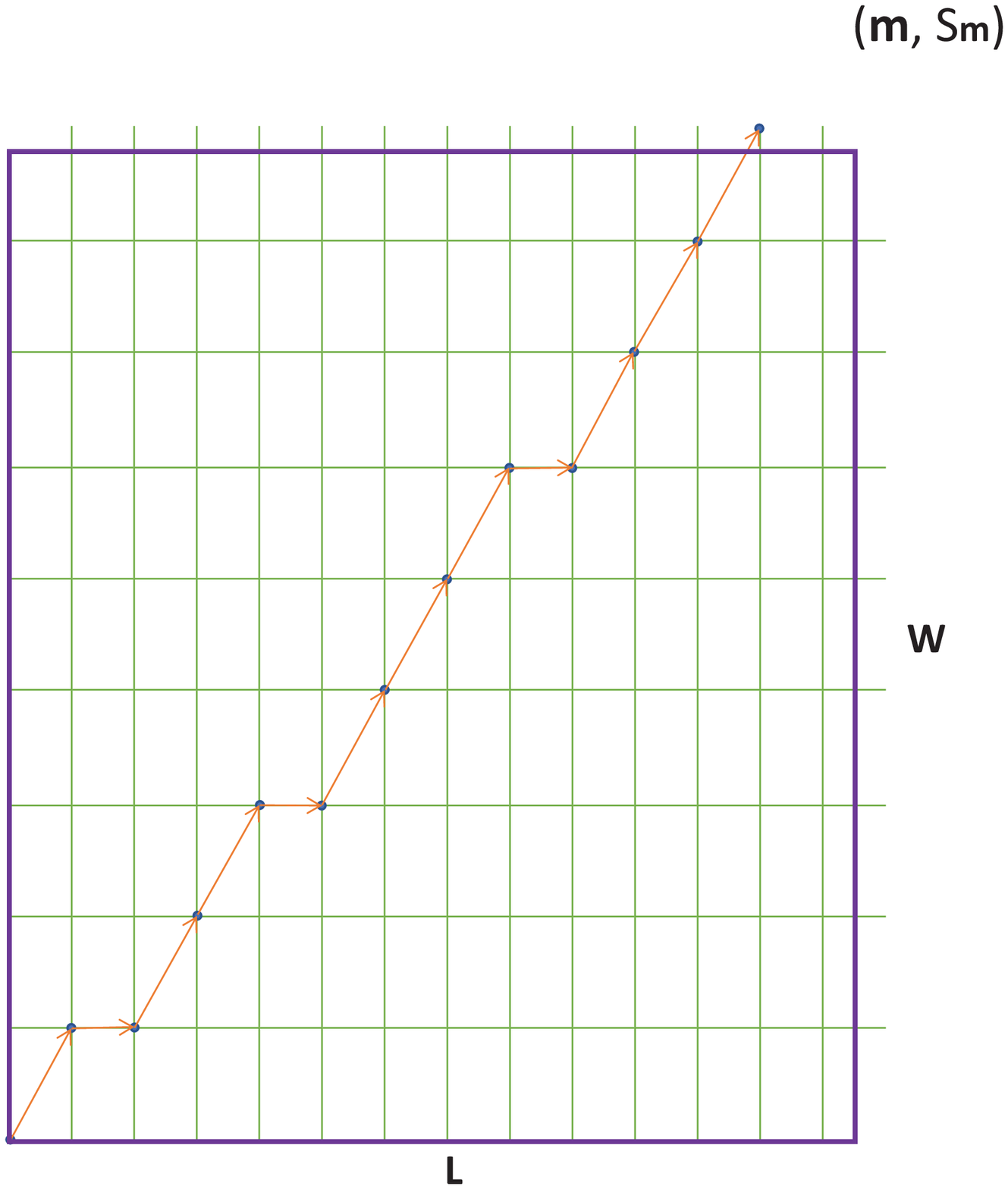, height=3in, width=5in }}
\caption{Rectangular Random Walk} \la{fig_Curve_A0}
\end{figure}

\bsk

Define the relative frequency \[ \wh{\bs{p}} = \f{ S_{\bs{m}} }{
\bs{m} }
\]
 as an estimator for the desired probability  $p$.
According to Theorem \ref{themmainB}, the mixed criterion
(\ref{mixedc}) is satisfied.  We call the above random walk
 as Rectangular Random Walk as the point with coordinate $(n, S_n)$
 is moving until it is out of the rectangular domain $\bb{B}$.

\section{Worst-Case Analysis}

In this section, we shall investigate the worst-case performance of
the method for probability estimation proposed in Section 4.
Clearly, the maximum number of simulations is no greater than $L +
1$. This is a good thing, since the user can plan the worst
requirement of computational resources.

Most importantly, the proposed method is remarkably efficient. Under
the assumption that the margin of relative error $\ba$ is small, we
can show that the worst-case improvement upon the Chernoff-Hoeffding
bound is about a quarter of the ratio of the relative error margin
to the absolute error margin, that is, \be \la{GainASN}
 \f{ \tx{ Number of simulations by
Chernoff-Hoeffding Bound } } { \tx{Maximum number of simulations of
the proposed method} } \ap \f{ 1 }{ 4 } \times \f{ \ba }{\al}. \ee
As a consequence of the stopping rule for the retangular random
walk, the required number of simulations of our method can be much
smaller than the worst case bound $L + 1$. Hence, the average
improvement of our method upon the Chernoff-Hoeffding bound can be
much greater than \[ \f{ 1 }{ 4 } \times \f{ \ba }{\al}.
\]
To justify the approximate formula (\ref{GainASN}), we employ Taylor
series expansion theory to investigate the ratio
\[
\f{ \tx{ Number of simulations by Chernoff-Hoeffding Bound } } {
\tx{Maximum number of simulations of the proposed method} }.
\]
Clearly,
\[
\f{ \tx{ Number of simulations by Chernoff-Hoeffding Bound } } {
\tx{Maximum number of simulations of the proposed method} } \ap \f{
\f{ \ln \f{2}{\de} }{ 2 \al^2 } }{ L  },
\]
where
\[
L = \f{\ba}{ (1 + \ba) \ln (1 + \ba) - \ba} \f{ \ln \f{2}{\de}
}{\al}.
\]
Hence,
\[
\f{ \tx{ Number of simulations by Chernoff-Hoeffding Bound } } {
\tx{Maximum number of simulations of the proposed method} } \ap
\f{1}{2} \f{\ba}{\al} \times \f{ (1 + \ba) \ln (1 + \ba) - \ba
}{\ba^2}.
\]
By Taylor series expansion theory, we have
\[
\ln ( 1 + x) \ap x - \f{x^2}{2}
\]
for small $|x|$.   Hence,
\[
\f{ (1 + \ba) \ln (1 + \ba) - \ba }{\ba^2} \ap \f{  (1 + \ba) ( \ba
- \f{\ba^2}{2} ) - \ba }{ \ba^2 } = \f{1}{2} (1 - \ba) \ap \f{1}{2},
\]
which implies (\ref{GainASN}).

The approximate formula (\ref{GainASN}) indicates that our proposed
method offers an extremely significant advancement in terms of
efficiency, without loss of rigorousness. The rationale is that in
practices, one can accept a relative margin in orders of magnitude
larger than an absolute margin. For example, to estimate the
probability $p$ of a critical failure, it is expected to estimate
$p$ within an absolute error of $10^{-6}$. The required number of
simulations calculated  with the Chernoff-Hoeffding bound is
$3,800,451,229,772$ for $\de = 10^{-3}$. Even with current super
computers, it is unthinkable to perform such an astronomical number
of simulations. To overcome such difficulty of computational
complexity, we use the mixed criterion proposed in Section 4.1.  We
consider a very mild relaxation by introducing a requirement of
relative error. According to the approximate formula
(\ref{GainASN}), our method can lead to a reduction of computation
by a factor of
\[ \f{1}{4} \times \f{10^{-2}}{10^{-6}} = 2,500. \]
Note that an estimate of one percent relative error is an extremely
stringent requirement of accuracy. More reduction of computation is
possible if the margin of relative error is increased.  For example,
if the margin of relative error is $0.1$, the computation can be
reduced by a factor of
\[ \f{1}{4} \times \f{10^{-1}}{10^{-6}} = 25,000. \]

\section{Conclusion}

We have developed a general theory of truncated inverse binomial
sampling.  The theory have been applied to estimate the probability
of event.  Worst-case analysis shows that the proposed method is
extremely efficient as compared to the widely used
Chernoff-Hoeffding bound, without scarifying the rigorousness of
error control in the estimation of probability of event.

\end{document}